\documentclass[11pt]{amsart}
\usepackage[T1]{fontenc}
\usepackage[utf8]{inputenc}
\usepackage{lmodern}
\usepackage{amsmath,amssymb,amsthm}
\usepackage{geometry}
\usepackage[hidelinks]{hyperref}
\usepackage{needspace}
\usepackage{etoolbox}
\usepackage{enumitem}
\usepackage{booktabs}
\setlist[enumerate,1]{label=\textup{(\alph*)}}

\makeatletter
\def\subsection{\@startsection{subsection}{2}%
  \z@{.5\linespacing\@plus.7\linespacing}{.25\linespacing}%
  {\normalfont\bfseries}}
\def\subsubsection{\@startsection{subsubsection}{3}%
  \z@{.5\linespacing\@plus.7\linespacing}{.2\linespacing}%
  {\normalfont\itshape}}
\makeatother

\pretocmd{\section}{\needspace{4\baselineskip}}{}{}%
\pretocmd{\subsection}{\needspace{3\baselineskip}}{}{}%

\newcommand{\Ord}{\mathrm{Ord}}

\newcommand{\ZFC}{\mathsf{ZFC}}

\newcommand{\PA}{\mathsf{PA}}

\newcommand{\NFUA}{\mathsf{NFUA}}
\newcommand{\ZF}{\mathsf{ZF}}

\theoremstyle{definition}
\newtheorem{definition}{Definition}[section]
\theoremstyle{plain}
\newtheorem{theorem}[definition]{Theorem}

\newtheorem{corollary}[definition]{Corollary}
\theoremstyle{remark}

\AtBeginEnvironment{theorem}{\setlist[enumerate]{itemsep=.4\baselineskip}}
\AtBeginEnvironment{proposition}{\setlist[enumerate]{itemsep=.4\baselineskip}}
\AtBeginEnvironment{lemma}{\setlist[enumerate]{itemsep=.4\baselineskip}}
\AtBeginEnvironment{corollary}{\setlist[enumerate]{itemsep=.4\baselineskip}}
\begin{document}

\title{What David Teaches Goliath: \\ An Inevitable Extension of ZFC}
\author{Ali Enayat}
\address{Department of Philosophy, Linguistics, and the Theory of Science\\
University of Gothenburg\\
Gothenburg, Sweden}
\email{ali.enayat@gu.se}
\date{}
\subjclass[2020]{Primary 03C62; Secondary 03E30, 03E55, 03E70}
\keywords{Levy scheme, Mahlo cardinals, reflection principles, models of set theory,
elementary end extensions, weak compactness of Ord, indiscernibles, NFUA}

\begin{abstract}
We examine the metamathematical significance of the Levy scheme $\Lambda$, which asserts, for each natural number $n$, the existence of an $n$-Mahlo cardinal $\kappa$ such that $V_\kappa\prec_{\Sigma_n}V$. Several seemingly unrelated results arrive at $\mathsf{ZFC}+\Lambda$: comparisons with the model theory of arithmetic and finite set theory; class theories in which $\Ord$ has the tree property; set theories with a proper class of indiscernibles; and interpretations of ordinary set theory in Cantorian parts of models of $\mathsf{NFUA}$, an extension of Quine--Jensen set theory $\mathsf{NFU}$.

Taken together, these results suggest that $\mathsf{ZFC}+\Lambda$ is an inevitable extension of $\mathsf{ZFC}$ when certain metamathematical phenomena familiar from the microcosm of hereditarily finite sets are transported to the macrocosm of all sets.

\end{abstract}

\maketitle
\pagestyle{plain}
\raggedbottom

\section{Introduction}

One of the recurrent themes in the foundations of set theory, going
back at least to G\"odel's reflections on new axioms for set theory
\cite{GodelContinuum}, is the search for natural extensions of
$\mathsf{ZFC}$. This paper highlights the role of Mahlo cardinals of finite order and corresponding reflection schemes in
the foundations of set theory by relying on purely metamathematical findings
to bring attention to the singular status of an axiom scheme $\Lambda $
(here dubbed the \textit{Levy scheme}) in the axiomatization of set
theory. The Levy scheme is a collection of first order sentences formulated in the
usual vocabulary $\{\in ,=\}$ of set theory that asserts, for each natural number $n$ (in the metatheory), the existence of $%
\Sigma _{n}$-reflective $n$-Mahlo cardinals. Since $ \Sigma _{n}$-reflective cardinals are increasingly referred to as $n$-correct cardinals, we will follow this practice in the rest of the paper.

In the elaborate hierarchy of large cardinals, inaccessible and Mahlo cardinals,
along with their transfinitely iterated generalizations ($\alpha $%
-inaccessible cardinals, and $\alpha $-Mahlo cardinals) belong to the
junior league of \textit{small} large cardinals, residing far below the
dominion of the so-called \textit{medium} large cardinals (such as weakly
compact cardinals and various flavors of ineffable and subtle cardinals).
Small and medium large cardinals are in turn dwarfed in comparison with 
\textit{large} large cardinals, cardinals whose existence contradict G\"{o}%
del's axiom of constructibility (such as measurable, supercompact, Woodin, and huge
cardinals). It is fair to say that large large cardinals have attracted the
lion's share of attention from set theorists during the last two or three
decades, often with spectacular results. However, as we shall argue, there
is still a great deal to be learnt about small large cardinals. Indeed, the
work of Harvey Friedman on Borel Diagonalization \cite{Friedman1981} exhibits the unexpected
impact of Mahlo cardinals of finite order on the tangible realm of real
analysis. Friedman has a large body of work devoted to the 
ramifications of large cardinals, especially of the Mahlo variety, for
``ordinary'' mathematical practice; for a representative sample, see \cite{Friedman2004}. The published results surveyed in this paper  complement Friedman's work by
bringing attention to new internal motives within the discourse of set theory and
model theory that indicate the  inevitability of $n$-Mahlo cardinals that are also $n$-correct, as codified in the Levy scheme. 

We will survey a range of metamathematical results 
concerning finite set theory when they are transported to ordinary set theory.  Given some brand of set or class theory, $T$, which includes the axiom of infinity, we use $T^{-\infty}$ to denote the "finitistic" variant of $T$ in which $\mathsf{Inf}$ is replaced by its negation and the axiom $\mathsf{TC}$, asserting the existence of transitive closures, is added.\footnote{The reason is technical, and is motivated by the indispensable role of  $\mathsf{TC}$ in the bi-interpretability of Peano arithmetic and finite set theory. For more detail, see \cite{KW} and \cite{ESV}.} The principal examples discussed here are shown below:

\[
\renewcommand{\arraystretch}{1.55}
\begin{array}{@{\quad}c@{\qquad}c@{\quad}}
\toprule
\text{\bfseries Finite set/class theory}
&
\text{\bfseries Infinite set/class theory}
\\
\midrule
\mathsf{ZFC}^{-\infty}
&
\mathsf{ZFC}+\Lambda
\\
\mathsf{GBC}^{-\infty}+\text{``}\Ord\text{ is weakly compact''}
&
\mathsf{GBC}+\text{``}\Ord\text{ is weakly compact''}
\\
\mathsf{ZFCI}_{<}^{-\infty}
&
\mathsf{ZFCI}_{<}
\\
\mathsf{NFUA}^{-\infty}
&
\mathsf{NFUA}
\\
\bottomrule
\end{array}
\]

\medskip

As we will see, the rows of the table above are related as follows on the basis of published results, with the finite $\mathsf{NFUA}$ characterization in Exhibit IV drawn from unpublished joint work with Solovay.
\begin{itemize}
\item Various well-known model theoretic results (the behavior of end extensions, existence of powerlike models, and well-behaved automorphisms) about completions of $\mathsf{ZFC}^{-\infty}$  carry over to a completion $T$ of $\ZFC$ iff $T\vdash \Lambda $; see Exhibit I (Section~\ref{sec:exhibit-I}).
    \item The first-order set-theoretic consequences of
$\mathsf{GBC}^{-\infty}~+~$``$\Ord$ is weakly compact'' are exactly those of
$\mathsf{ZFC}^{-\infty}$, but the first-order set-theoretic consequences of
$\mathsf{GBC}~+~$``$\Ord$ is weakly compact'' are exactly those of
$\mathsf{ZFC}+\Lambda$; see Exhibit II (Section~\ref{sec:exhibit-II}).
\item $\mathsf{ZFCI}^{-\infty}_{<}$ ($\mathsf{ZFC}^{-\infty}$ with a class of indiscernibles) is conservative over $\mathsf{ZFC}^{-\infty}$,
whereas the pure set-theoretic consequences of 
$\mathsf{ZFCI}_{<}$ ($\mathsf{ZFC}$ with a class of indiscernibles) are precisely those of $\mathsf{ZFC}+\Lambda$; see Exhibit III (Section~\ref{sec:exhibit-III}). 
\item The complete theories of the Cantorian parts of models of
$\mathsf{NFUA}^{-\infty}$ are precisely the completions of
$\mathsf{ZFC}^{-\infty}$, whereas the complete theories of the
Cantorian parts of models of $\mathsf{NFUA}$ are exactly the
completions of $\mathsf{ZFC}+\Lambda$; see Exhibit IV (Section~\ref{sec:exhibit-IV}).
\end{itemize}

Although the four exhibits concern seemingly disparate aspects of set theory, their proofs reveal some noteworthy connections. In particular, the automorphism and indiscernibility constructions both make use of strong cuts, a notion originating in the model theory of arithmetic, in obtaining class expansions in which $\Ord$ is weakly compact. The characterization of the Cantorian parts of models of $\mathsf{NFUA}$, in turn, relies on the automorphism constructions. These connections will be explained in the respective exhibits.

There are two important points to make clear at this point:

\begin{enumerate}
    \item The axiom of choice is provable in $\ZF^{-\infty}$; indeed the statement $\mathsf{V=L}$ is provable in $\ZF^{-\infty}$. Consequently, the deductive closure of  $\ZF^{-\infty}$ coincides with the deductive closure of $\ZFC^{-\infty}$. Therefore, for the purposes of this paper, they are interchangeable. 
    \item The above also makes the distinction between $\mathsf{GB}^{-\infty}$ and $\mathsf{GBC}^{-\infty}$ immaterial here. Combined with the fact that the proof of König's tree lemma can be readily formalized in $\mathsf{GB}^{-\infty}$, this shows that  $\mathsf{GB}^{-\infty}$, $\mathsf{GBC}^{-\infty}$, and $\mathsf{GBC}^{-\infty} + \text{``}\Ord\text{ is weakly compact''}$ have the same deductive closure. Therefore, for the purposes of this paper, they are interchangeable.
\end{enumerate}

There is a suggestive way of understanding why the comparison between finite set theory and full set theory leads toward large-cardinal principles such as $\Lambda$. As observed in
\cite{E-3}, the negation of Infinity may itself be regarded, in a qualified
sense, as a large-cardinal principle: rather than positing the existence
of a large set---indeed, it denies the existence of any infinite set---it
attributes large-cardinal-like properties to the universe as a whole.
Power Set and Replacement, by contrast, give the class of ordinals a
strong-inaccessibility character. This helps explain why some
metamathematical phenomena familiar from arithmetic and finite set theory
reappear in infinite set theory only under additional large-cardinal
hypotheses; in the above comparisons, the additional strength
is precisely that codified by $\Lambda$.

\medskip

\noindent\textit{A personal note.}
In a sense, I have been writing this paper since 2009, when I gave a talk \cite{E-2009}
announcing the key results surveyed here, some, but not all, of which had by
then appeared in print. One reason it took me so long to write the paper was
the need to work out a satisfactory account of various technical details of the theory $\ZFC^{\pm\infty}$ with a distinguished proper class of indiscernibles, which I managed to finally accomplish, not too long ago, in \cite{E-ZFI, E-PAI}. In preparing the
present paper I have committed what might be called a considerable amount of
self-plagiarism, freely drawing on my earlier writings, with the assistance of
the edu version of ChatGPT (courtesy of University of Gothenburg), which came to my rescue in organizing this material and in editing and typesetting the \LaTeX{} manuscript.  

\section{Preliminaries}

Here we collect the basic definitions, notations, conventions, and results
that are relevant for the subsequent sections.  
\medskip

\subsection{Theories and models}

\begin{itemize}
\item  Given a language $\mathcal{L}\supseteq \{\in
,=\}, $ $\Sigma _{n}(\mathcal{L}),$ $\Pi _{n}(\mathcal{L})$, and $\Delta
_{n}(\mathcal{L})$ are sets of first order formulas built up from $\mathcal{L%
}$ that belong, respectively, to the $\Sigma _{n}$-level, $\Pi _{n}$-level,
and $\Delta _{n}$-level of the well-known Levy hierarchy (see \cite{Jec,Kan}).\ We use $%
\Sigma _{\infty }(\mathcal{L})$ for $\bigcup\limits_{n<\omega}\Sigma
_{n}(\mathcal{L}).$

\item  For a class $\Gamma $ of formulas, $\Gamma $-\textit{Separation }%
consists of the universal closure of formulas of the form 
\[
\forall a\exists b\forall x(x\in b\longleftrightarrow x\in a\wedge \varphi
(x,v)), 
\]
for all formulas $\varphi \in \Gamma $ in which $b$ does not occur free
(note: $v$ is a parameter).

\item  For a class $\Gamma $ of formulas, $\Gamma $-\textit{Replacement }%
consists of the universal closure of formulas of the form 
\[
[\forall x\in a\,\exists!y\,\varphi(x,y,v)]\rightarrow
[\exists b\,\forall y\,(y\in b\leftrightarrow
\exists x\in a\,\varphi(x,y,v))],
\]
for all formulas $\varphi \in \Gamma $ in which $b$ does not occur free.

\item  $\mathsf{ZFC}(\mathcal{L})$ is $\mathsf{ZFC}$ augmented by $\Sigma _{\infty }(\mathcal{L}%
) $-Replacement.\ 

\item  Suppose $\mathcal M=(M,E)$\ is a model of $\mathsf{ZFC}$, where $E=$ $\in ^{%
\mathcal M}$.

\textup{(a)} $\Ord^{\mathcal M}$ denotes the ordered set of \ ``ordinals''
of $\mathcal M$.

\textup{(b)} If $\mathcal M\subseteq $ $\mathcal N\ =(N,F)$, then $\mathcal N$ 
\textit{end extends} $\mathcal M$, written $\mathcal M\subseteq _{e}%
\mathcal N$, if $\mathcal N$ fixes every element of $\mathcal M$,
i.e., for every $c\in M,$ $\{x\in M:xEc\}=\{x\in N:xFc\}.$

If, in addition, $\mathcal M\prec\mathcal N$, then $\mathcal N$ is an
\textit{elementary end extension} of $\mathcal M$. An end extension
$\mathcal N$ of $\mathcal M$ is a \textit{rank extension} if every member of
$N\setminus M$ has rank, in $\mathcal N$, above the rank of every member of
$M$. An elementary extension $\mathcal N$ of $\mathcal M$ is \textit{minimal}
if there is no proper intermediate elementary submodel $\mathcal K$ with
$\mathcal M\prec\mathcal K\prec\mathcal N$. An extension $\mathcal N$ of
$\mathcal M$ is \textit{conservative} if, for every subset $A$ of $N$
definable in $\mathcal N$ (with parameters), $A\cap M$ is definable in
$\mathcal M$ (with parameters).

\item  $L_{\alpha }$ is the $\alpha $-th level of G\"{o}del's class of
constructible sets $L$. If $\mathcal N$ is also a model
of $\mathsf{ZFC}$ with $\mathcal M\subseteq _{e}\mathcal N$ and $\Ord^{%
\mathcal M}=\Ord^{\mathcal N}$, then $L^{\mathcal M}=%
L^{\mathcal N}$.

\item  $\rho $ is the usual ordinal valued rank function on sets and $%
V_{\alpha }$ is the $\alpha $-th level of the von Neumann cumulative
hierarchy, so $V_{\alpha }=\{x:\rho (x)<\alpha \}.$

\item  Our blanket assumption throughout the paper is that $\mathsf{ZFC}$\textit{\ is
a consistent first order theory}.\ Also, a \textit{completion} of a
consistent theory $T$ refers to a \textit{consistent} completion of $%
T. $
\end{itemize}

\subsection{Mahlo cardinals and reflection}

\begin{itemize}
\item  A cardinal $\kappa $ is a \textit{strongly inaccessible cardinal} if $%
\kappa $ is regular and strong limit. As observed by Kuratowski in the
1930's, the strong inaccessibility of $\kappa $ implies that $V_{\kappa }$
is a model of $\mathsf{ZFC}.$ Therefore $\mathsf{ZFC}$ cannot prove the existence of a strongly
inaccessible cardinal. Furthermore, G\"{o}del's incompleteness theorem can
be used to show that the statement Con$(\mathsf{ZFC}+$ ``there is an inaccessible
cardinal'') is not provable within $\mathsf{ZFC}$, assuming $\mathsf{ZFC}+$ Con$(\mathsf{ZFC})$ is
consistent\footnote{%
For a recursively axiomatizable theory $T$, Con($T$)\ is the \textit{%
arithmetical statement} expressing ``$T$ is consistent''.}.

\item  A Mahlo cardinal $\kappa $ is a strongly inaccessible cardinal such that the
set of regular cardinals below $\kappa$ is \textit{stationary}
in $\kappa$ (i.e., it meets every closed and unbounded subset of $%
\kappa$). The $\alpha$-Mahlo cardinals, for an ordinal $\alpha$,
are defined recursively as follows:

$\kappa $ is 0-Mahlo if $\kappa $ is strongly inaccessible; and

$\kappa $ is $(\alpha +1)$-Mahlo if $\kappa $ is a strongly inaccessible
cardinal and $\{\lambda <\kappa :\lambda $ is $\alpha $-Mahlo$\}$ is
stationary in $\kappa .$

For a limit ordinal $\alpha >0,$ $\kappa $ is $\alpha $-Mahlo if for all $%
\beta <\alpha ,$ $\kappa $ is $\beta $-Mahlo.

\item For $n<\omega$, an ordinal $\kappa$ is \emph{$n$-correct} (or \emph{$\Sigma_n$-correct}) if
\[
V_\kappa\prec_{\Sigma_n}V.
\]
Thus the Levy scheme considered below synchronizes two finite hierarchies:
finite Mahloness and finite correctness.$\bigskip $
\end{itemize}

\subsection{Reflection theorems}

\begin{itemize}
\item  Levy showed that, in contrast with Tarski's undefinability of
Truth, $\Sigma _{n}$-truth is $\Sigma _{n}$-definable for $n\geq 1$ within $%
\mathsf{ZFC}$ \cite[Section 14]{Jec}. Therefore for each natural number $n$ there is a unary
first order formula in the language of set theory with the free variable $%
\alpha $, denoted ``$V_{\alpha}\prec_{\Sigma_n}V$'',
that expresses the statement ``for all $\Sigma_n$-formulae
$\varphi(v_0,\ldots,v_k)$ and all $a_0,\ldots,a_k\in V_\alpha$,
$\varphi(a_0,\ldots,a_k)$ holds if and only if
$\varphi^{V_\alpha}(a_0,\ldots,a_k)$ holds'', where
$\varphi^{V_\alpha}$ denotes the relativization of $\varphi$
to $V_\alpha$.

\item For each formula $\varphi(x,\vec y)$ in the language of set theory,
we write $\operatorname{Club}_{\varphi}(\vec a)$ for the first-order
formula asserting that $\{\alpha\in\Ord:\varphi(\alpha,\vec a)\}$
is a closed unbounded class of ordinals.\medskip
\end{itemize}

The following theorem is Levy's fine-tuned version of Montague's
reflection theorem. 
This version follows from the usual version in \cite{Jec,Kan} by simply observing that the set of $\Sigma _{n}$-formulae is closed under subformulae
and therefore, by Tarski's test, $\prec_{\Sigma_n}$ is preserved in limits.\medskip

\begin{theorem}[Montague--Levy]
For each $n < \omega$,
\[
\mathsf{ZFC} \vdash \text{``}\{\alpha \in \operatorname{Ord} \mid V_{\alpha} \prec_{\Sigma_n} V\} \text{ is a club''}.
\]
\end{theorem}

By the following classical theorem of Levy, one obtains \textit{%
full-reflection} at Mahlo cardinals of finite order:

\begin{theorem}[Levy]\textit{For each} $n<\omega$ \textit{the following statement is provable
within} $\mathsf{ZFC}$:\smallskip

\begin{center}
$(\kappa $ \textit{is} $(n+1)$-\textit{Mahlo}) $\rightarrow $ $\exists
\alpha <\kappa $ $(\alpha $ is $n$-Mahlo and $V_{\alpha }\prec V_{\kappa
}).\medskip $
\end{center}
\end{theorem}

Note that the results in this section are valid for $\mathsf{ZFC}(%
\mathcal{L})$, for a finite relational $\mathcal{L}\supseteq \{=,\in \}.$

\section{Basic facts about \texorpdfstring{$\mathsf{ZFC}+\Lambda$}{ZFC + Lambda}}

In this section we officially introduce the Levy scheme $\Lambda $ and
discuss fundamental metamathematical aspects of $\mathsf{ZFC}+\Lambda $. We first examine alternative axiomatizations of $\Lambda$ over $\mathsf{ZFC}$, and discuss its relationship with the weaker scheme $\Lambda^{-}$, defined below.

We then establish the robustness of $\Lambda$ by
showing that $\Lambda $ is invariant under two fundamental transformations
in set theory: the move \textit{inward} to the constructible universe,
and the move \textit{outward} to a generic extension.\medskip

\subsection{\texorpdfstring{$\Lambda$}{Lambda} and its equivalents}

For a unary predicate $X$ and each $n<\omega$, let $\lambda_n(X)$ assert
that there is an $n$-Mahlo cardinal $\kappa$ such that
\[
 (V_\kappa,\in,X\cap V_\kappa)\prec_{\Sigma_n(X)}(V,\in,X).
\]
We write
\[
\Lambda(X):=\{\lambda_n(X):n<\omega\},
\]
and use $\Lambda$ for the scheme without $X$. Also let
\[
\Lambda^{-}:=\{\text{``there is an $n$-Mahlo cardinal''}:n<\omega\}.
\] 

The passage from $\Lambda^{-}$ to $\Lambda$, by requiring the $n$-Mahlo cardinals in question to be also $n$-correct, turns a sequence of assertions about large cardinals below $\Ord$ into a reflection principle naturally attributable to $\Ord$ itself. It is in this qualified sense that $\Lambda$ can be read schematically as saying that ``$\Ord$ is $\omega$-Mahlo''.

\medskip

Thanks to Levy's pioneering work, the following theorem is now straightforward to prove.

\begin{theorem}[Levy]
Over $\mathsf{ZFC}$, the following schemes are deductively equivalent:
\begin{enumerate}
\item[\textup{(a)}] $\Lambda$;

\item[\textup{(b)}] The scheme
$\{\text{``there are arbitrarily large $n$-Mahlo cardinals
$\kappa$ such that $V_\kappa\prec_{\Sigma_n}V$''}:n<\omega\}$;

\item[\textup{(c)}] The scheme
$\{\text{``there are arbitrarily large $m$-Mahlo cardinals
$\kappa$ such that $V_\kappa\prec_{\Sigma_n}V$''}:m,n<\omega\}$;

\item[\textup{(d)}] The scheme consisting, for each $n<\omega$
and each formula $\varphi(x,\vec y)$ in the language
of set theory, of the sentence
\[
\forall\vec a\,
\left(
\operatorname{Club}_{\varphi}(\vec a)
\longrightarrow
\exists\kappa\,
\bigl(
\varphi(\kappa,\vec a)\land
\kappa\text{ is }n\text{-Mahlo}
\bigr)
\right).
\]
\end{enumerate}
\end{theorem}

The next result is \cite[Theorem 2.4.11(c),(d)]{E-ZFI}.
\begin{theorem}
$\mathsf{ZFC}+\Lambda^{-}$ and $\mathsf{ZFC}+\Lambda$ are mutually interpretable.\footnote{However, the two theories are not bi-interpretable. This follows from a more general result asserting that any two bi-interpretable extensions of $\ZF$ (in the same language as $\ZF$) have the same deductive closure. See \cite{E-Visser,E-Visser-Corr}.}
Assuming the consistency of $\mathsf{ZFC}+$ ``there is an $\omega$-Mahlo cardinal'',
$\mathsf{ZFC}+\Lambda^{-}$ does not prove $\Lambda$.
\end{theorem}

\medskip

\subsection{The robustness of \texorpdfstring{$\Lambda$}{Lambda}}

\medskip

The following theorem has been known to the cognoscenti since the 1960s; a published formulation appears as \cite[Theorem 2.4.11(e)]{E-ZFI}.
\begin{theorem}[Relative constructibility robustness]
If $\mathcal M\models\mathsf{ZFC}+\Lambda$, then for every
$c\in M$, $L(c)^{\mathcal M}\models\Lambda$.
\end{theorem}

The following result is part of the folklore of the subject; for a proof see \cite[Theorem 2.4.11(f)]{E-ZFI}. It is morally due to Levy and Solovay \cite{LevySolovay}.
\begin{theorem}[Set-forcing preservation]
Suppose $\mathcal M\models \mathsf{ZFC}+\Lambda$, $\mathbb P\in M$, and $G$
 is $\mathbb P$-generic over $\mathcal M$. Then $\mathcal M[G]\models\Lambda$.
\end{theorem}
In contrast, there are class notions of forcing that destroy $\Lambda$ while preserving $\mathsf{ZFC}$. The underlying methods for altering Mahlo properties go back to a construction of Jensen \cite{JensenSolovay}, subsequently generalized by Boos \cite{Bo}.

\begin{corollary}
Suppose $\mathsf{ZFC}+\Lambda$ is consistent. If $\psi$ is a sentence in the language of set theory satisfying either

\begin{enumerate}
\item $\mathsf{ZFC}\vdash$ ``$\psi$ holds in $L$''; or

\item $\mathsf{ZFC}\vdash$ ``for some poset $\mathbb P$, $1_{\mathbb P}\Vdash\psi$''.
\end{enumerate}

Then there is a model of $\mathsf{ZFC}+\Lambda+\psi$.
\end{corollary}

\noindent Note that part (b) of the above corollary shows that if $\mathsf{ZFC}+\Lambda$ is consistent, then the continuum hypothesis is both consistent with and independent of $\mathsf{ZFC}+\Lambda$.

\section{Exhibit I: comparative model theory}
\label{sec:exhibit-I}

\subsection{Elementary end extensions }

It is well known \cite{KW} that $\PA$ (Peano Arithmetic) and $\mathsf{ZF}^{-\infty}$ are bi-interpretable. This allows one to translate many metamathematical results about $\PA$ into corresponding results concerning $\mathsf{ZF}^{-\infty}$, including the ones considered below. The underlying coding goes back to Ackermann \cite{Ackermann}; see also Mycielski \cite{Mycielski} and Kaye--Wong \cite{KW} for the development of the correspondence. The following result was formulated for $\PA$, but thanks to the aforementioned bi-interpretability result, it can equally well be presented as below.

\begin{theorem}
\textnormal{(MacDowell--Specker \cite{MacDowell-Specker}, Gaifman
\cite{G}, Phillips \cite{Phillips}).}
Every model $\mathcal M$ of $\mathsf{ZF}^{-\infty}$ has a proper
elementary end extension $\mathcal N$.  Moreover, $\mathcal N$ can be
required to be both minimal and conservative over $\mathcal M$.
\end{theorem}

The above result naturally prompts the question of whether the result holds if $\ZF^{-\infty}$ is replaced by $\ZF$. The following theorem provides a partial answer to this question (in a strong form).

\begin{theorem}[Keisler--Morley \cite{Keisler-Morley}]
Suppose $\mathcal M$ is a model of $\mathsf{ZF}$ whose ordinal order has
cofinality $\aleph_0$, and let $\langle X,\lhd\rangle$ be a linear order
of arbitrary cardinality.  There is an elementary end extension
$\mathcal N$ of $\mathcal M$ whose ordinal order contains an isomorphic
copy of $\langle X,\lhd\rangle$ as indiscernibles above the ordinals of
$\mathcal M$.
\end{theorem}

The countable-cofinality hypothesis marks an essential limitation: unlike
the finite situation, not every model of $\mathsf{ZFC}$ has a proper
elementary end extension. Keisler and Silver \cite{Keisler-Silver} showed that, if
$\kappa$ is the first inaccessible cardinal, then $V_\kappa$ has no proper
elementary end extension.  Subsequent results of Kaufmann and Enayat
\cite{Kau,E-1,E-3} showed that every consistent completion of
$\mathsf{ZFC}$ has a model of power $\aleph_1$ with no proper elementary end extension. The strongest result in this direction is as follows.

\begin{theorem}
\textnormal{(\cite[Theorem C]{E-Deadends}).}
Every consistent extension $T$ of $\mathsf{ZF}$ has a model $\mathcal M$
of cardinality $\aleph_1$ such that $\mathcal M$ has no proper end
extension to a model of $\mathsf{ZF}$.
\end{theorem}

Thus the full MacDowell--Specker
analogue fails dramatically for set theory. The appropriate positive analogue therefore has to be formulated more
carefully.  Instead of requiring every model of a completion $T$ itself
to have such an extension, one passes to a suitable \emph{countable expansion
of the language} and asks for a consistent expanded theory $\overline T$ all of whose
models have an elementary end extension.  At precisely this point the provability of the Levy scheme in $T$ emerges as a necessary and sufficient condition.

\begin{theorem}\label{Lambda-eee}
\textnormal{(\cite{Kau}, \cite {E-5}).}
Let $T$ be a consistent completion of $\mathsf{ZFC}$. The following are
equivalent:
\begin{enumerate}
\item[\textup{(a)}] $T\vdash\Lambda$.
\item[\textup{(b)}] There is a consistent theory $\overline T$ in a countable
expansion $\mathcal L$ of the language of set theory, containing a global
well-order, such that
\[
T+\mathsf{ZFC}(\mathcal L)\subseteq\overline T,
\]
and every model of $\overline T$ has a proper \emph{elementary
end extension}.
\item[\textup{(c)}] There is a consistent theory $\overline T$ in a countable
expansion $\mathcal L$ of the language of set theory, containing a global
well-order, such that
\[
T+\mathsf{ZFC}(\mathcal L)\subseteq\overline T,
\]
and every model of $\overline T$ has a proper \emph{conservative elementary
end extension}.
\end{enumerate}
\end{theorem}

Thus $\Lambda$ exactly marks the
completions for which this uniform expanded-language analogue of the
MacDowell--Specker phenomenon is available. See \cite[Section 2]{E-5} for the
construction and the converse argument.

\medskip

\subsection{Power-like models}

Recall that a model $\mathcal M=(M,E,\ldots)$ is called
$\theta$-\emph{like}, for an uncountable cardinal $\theta$, if
$|M|=\theta$ and every $E$-extension
$\{b\in M:bEa\}$, for $a\in M$, has cardinality less than
$\theta$. Thus the model itself has external size $\theta$, while every
set represented in the model has externally fewer than $\theta$ many
members. The terminology is meant to suggest the familiar situation of
$V_\theta$ when it is a model of set theory. Power-like models also have
a natural connection with two-cardinal transfer principles; see
\cite{E-4} for the surrounding results.

The point of departure of this section is the following result, which readily follows by putting the MacDowell--Specker theorem together with the bi-interpretation between $\PA$ and $\ZF^{-\infty}$.

\begin{theorem} Suppose $T$ is a completion of $\ZF^{-\infty}$. For every uncountable cardinal $\theta$, $T$ has a $\theta$-like model.
    
\end{theorem}

The contrast between $\aleph_1$ and higher cardinals is important here.
The Keisler--Morley theorem above implies that every consistent extension of
$\mathsf{ZF}$ has an $\aleph_1$-like model. At higher cardinals the
situation is subtler. For example, Chang's two-cardinal theorem implies,
in the presence of the continuum hypothesis, that every consistent
extension of $\mathsf{ZF}$ has an $\aleph_2$-like model, although CH is
not necessary for this conclusion \cite{E-4}. On the other hand, when
the cardinal at which a power-like model is sought has the tree
property, the existence of such a model imposes a strong restriction on
its complete theory.

\begin{theorem}
\textnormal{(\cite[Theorem A]{E-4}).}\label{Theorem-power-like}
Suppose $\theta$ is an uncountable cardinal with the tree property, and
$T$ is a consistent completion of $\mathsf{ZFC}$. The following are
equivalent:
\begin{enumerate}
\item[\textup{(a)}] $T$ has a $\theta$-like model.
\item[\textup{(b)}] $T$ contains the Levy scheme $\Lambda$.
\item[\textup{(c)}] $T$ has a $\kappa$-like model for every uncountable cardinal
$\kappa$.
\end{enumerate}
\end{theorem}

The proof of Theorem \ref{Theorem-power-like} in \cite{E-4} uses a family of ranked trees which, in the
terminology of that paper, ``try to elongate the universe.'' For a
$\theta$-like model, these trees are externally $\theta$-like, so the tree
property of $\theta$ supplies branches from which one constructs
$\Sigma_n$-elementary end extensions, for each $n<\omega$. By passing to
suitable extensions having a least new ordinal $\gamma$, one shows by
induction that $\gamma$ is $n$-Mahlo. The elementarity of the extensions
then yields the corresponding instances of $\Lambda$ in the original model.
A later variant of the tree construction appears in \cite{EH-Ord}.

Thus, at a cardinal with the tree property, the apparently semantic
condition that a completion possess a power-like model is equivalent to the
syntactic condition that it contain $\Lambda$. 

The role of the tree property also explains why $\aleph_2$ is a
particularly revealing test case. The mere existence of an
$\aleph_2$-like model does not by itself characterize $\Lambda$; what is
striking is what happens when one asks uniformly that all such models,
including finite relational expansions, satisfy the corresponding
Levy scheme.

\begin{theorem}
\textnormal{(\cite[Theorem B]{E-4}).}
The following assertions are equiconsistent over $\mathsf{ZFC}$:
\begin{enumerate}
\item[\textup{(a)}] There is an $\omega$-Mahlo cardinal.
\item[\textup{(b)}] For every finite language $\mathcal L$ expanding the language
of set theory, every $\aleph_2$-like model of
$\mathsf{ZFC}(\mathcal L)$ satisfies the corresponding expanded-language
scheme $\Lambda(\mathcal L)$.
\end{enumerate}
\end{theorem}

The two theorems therefore address different levels of the phenomenon.
The first characterizes those completions of $\mathsf{ZFC}$ having a
power-like model at a cardinal with the tree property; the second
measures the consistency strength of a uniform assertion about
\emph{all} $\aleph_2$-like models. The qualification concerning finite
expanded languages is essential.

\subsection{Automorphisms}

Gaifman refined the MacDowell-Specker method \cite{MacDowell-Specker} of building elementary end
extensions by introducing the powerful machinery of \textit{minimal end
extension types}, which can be used to produce\ a variety of models of $\mathsf{PA}$
with special properties. In particular, one can use minimal end extension
types to establish that \textit{every model }$\mathcal M=(M,+,\cdot )$%
\textit{\ of }$\mathsf{PA}$ \textit{has a proper elementary end extension} $%
\mathcal N$ \textit{that possesses an automorphism} $j$ \textit{whose
fixed point set is precisely }$M$\textit{\ }(this is
a corollary of \cite[Theorems 4.9--4.11]{G}). Therefore, in light of the bi-interpretability of $\PA$ and $\ZF^{-\infty}$ we have: \medskip

\begin{theorem}\label{thm:finite-automorphism}
Every consistent completion of $\ZF^{-\infty}$ has a model $\mathcal N$ that possesses a nontrivial
automorphism whose fixed point set is a proper initial segment of $\mathcal N$.
\end{theorem}

For the set-theoretic analogue, we use the notation of \cite{E-5}.
For a language $\mathcal L$, $\mathsf{EST}(\mathcal L)$ is Elementary Set
Theory in that language: Extensionality, Foundation, Pairing, Union,
Infinity, Choice, and $\Delta_0(\mathcal L)$-Separation. The axiom
$\mathsf{GW}$ says that the distinguished relation globally well-orders the
universe, and $\mathsf{GW}^{*}$ is the set-like, rank-respecting
strengthening used in \cite{E-5}. We refer to that paper for the precise
expanded-language axioms.

A natural question is whether there is an analogue of Theorem~\ref{thm:finite-automorphism} for $\mathsf{ZFC}$. This question has a positive answer \textit{precisely} for consistent
completions of the theory $\mathsf{ZFC}+\Lambda $. In other words:

\begin{theorem}
\textnormal{(\cite[Theorem A]{E-5}).}
Suppose $T$ is a consistent completion of $\mathsf{ZFC}+\Lambda$. There is a model $\mathcal M$
of $T+\mathsf{ZFC}(\vartriangleleft)+\mathsf{GW}$ such that $\mathcal M$
has a proper elementary end extension $\mathcal N$
possessing an automorphism whose fixed point set is precisely $M$.
\end{theorem}

\begin{theorem}
\textnormal{(\cite[Theorem B]{E-5}).}
Let $\mathcal N=(N,\in,\vartriangleleft)$ satisfy
$\mathsf{EST}(\{\in,\vartriangleleft\})+\mathsf{GW}^*$, and let $j$ be an
automorphism of $\mathcal N$ whose fixed-point set $M$
is a proper $\vartriangleleft$-initial segment of $N$.
There is a family $\mathcal A$ of classes over the induced
structure $\mathcal M=(M,\in)$ such that
\[
 (\mathcal M,\mathcal A)\models
 \mathsf{GBC}+\text{``$\Ord$ is weakly compact''}.
\]
Consequently, $\mathcal M\models \mathsf{ZFC}+\Lambda$.
\end{theorem}

A key ingredient in the converse direction is the notion of a
\emph{strong cut}, introduced by Kirby and Paris in the model theory of
arithmetic. In its set-theoretic analogue, a cut $M$ of an ordered model
$N$ is strong if, for every function $f$ coded in $N$ whose domain
includes $M$, there is an element $s\in N$ such that, for every $m\in M$,
\[
f(m)\notin M\quad\Longleftrightarrow\quad s<f(m).
\]
The proof of \cite[Theorem B]{E-5} shows that the fixed-point
structure of the automorphism is such a strong cut. This property is
used to show that the classes obtained by taking traces of sets in the
ambient model form a model of $\mathsf{GBC}$ over the fixed-point
structure, in which $\Ord$ is weakly compact. The Levy scheme then
follows by Exhibit II.
\medskip

\section{Exhibit II: G\"odel--Bernays theory of classes}
\label{sec:exhibit-II}

We now turn to a class-theoretic manifestation of the same phenomenon.  We
use $\mathsf{GBC}$ for G\"odel--Bernays class theory with global choice.  An
$\Ord$-tree is a class tree of height $\Ord$ whose levels are sets.  The
assertion
\[
             \text{``$\Ord$ is weakly compact''}
\]
means that every $\Ord$-tree has a cofinal class branch.  Since the branch
must belong to the class part of the model, this assertion is sensitive to
which classes accompany the sets; two models with the same first-order set
part can behave differently in this respect.

Over $\mathsf{GBC}$, the tree property for $\Ord$ is equivalent to the standard
class analogues of several familiar faces of weak compactness, including
infinitary compactness, partition properties, a Keisler-style extension
property, and $\Pi^1_1$-reflection; see \cite[Theorem 2.4.9]{E-ZFI} for the
precise formulations. At the same time one
must not transfer the cardinal picture too literally.  Within $\mathsf{ZFC}$, a weakly compact
cardinal is highly Mahlo, whereas the first-order set-theoretic strength of
``$\Ord$ is weakly compact'' turns out to stop exactly at the Levy scheme.
This is the point of the next theorem.

\begin{theorem}
\textnormal{(\cite[Theorem 2.1]{E-5}).}
\begin{enumerate}
\item[\textup{(a)}] If $(\mathcal M,\mathcal A)\models
$ $\mathsf{GBC}~+$ ``$\Ord$ is weakly compact'', then
$\mathcal M\models\mathsf{ZFC}+\Lambda$.
\item[\textup{(b)}] Every consistent completion $T$ of $\mathsf{ZFC}+\Lambda$ has a
countable model $\mathcal M$ admitting an expansion
$(\mathcal M,\mathcal A)\models\mathsf{GBC}+$ ``$\Ord$ is weakly compact''.
\end{enumerate}
\end{theorem}

The proof of part (a) uses the ranked-tree constructions encountered in Exhibit I. The tree property for $\Ord$ provides, at each finite level of the Levy hierarchy, suitable elementary end extensions with a least new ordinal. An induction on $n$ then shows that the least new ordinal is $n$-Mahlo, yielding $\Lambda$ in the original model. For part (b), the proof uses the Schmerl--Shelah method of constructing a theory in an expanded language whose models have conservative elementary end extensions. Taking the classes definable in the expanded structure yields the required $\mathsf{GBC}$ expansion.

\begin{corollary}
\textnormal{(\cite[Corollary 2.1.1]{E-5}; see also
\cite[Theorem 2.4.12]{E-ZFI}).}
For every sentence $\varphi$ of the pure language of set theory,
\[
\mathsf{GBC}+\text{``$\Ord$ is weakly compact''}\vdash\varphi
\quad\Longleftrightarrow\quad
\mathsf{ZFC}+\Lambda\vdash\varphi.
\]
\end{corollary}

There is also a useful warning about the slogan that $\Lambda$ says,
schematically, that ``$\Ord$ is $\omega$-Mahlo''.  \cite{E-5}
shows that although $\mathsf{GBC}+$ ``$\Ord$ is weakly compact'' proves each
standard finite level of Mahloness for $\Ord$, it does not prove the single
internal assertion that $\Ord$ is $\omega$-Mahlo.  The distinction is the
same scheme-versus-single-sentence distinction that is important throughout
this paper.

Finally, Enayat and Hamkins \cite{EH-Ord} showed that in every model of
$\mathsf{ZFC}$ the \emph{definable} $\Ord$-tree property fails: there is a
definable $\Ord$-tree with no definable cofinal branch.  This does not
contradict the theorem above.  In a $\mathsf{GBC}$ expansion witnessing weak
compactness, the available classes need not be restricted to the definable
classes of the set-theoretic reduct.  The result instead underscores why the
class part of the model is essential.

\section{Exhibit III: proper class of indiscernibles}
\label{sec:exhibit-III}
Indiscernibles have played an important role in model theory, starting from the pioneering work of Ehrenfeucht and Mostowski. In the set-theoretic realm, indiscernibles made their debut in the context of large cardinals and the study of $0^\#$. We now examine what happens when a proper class of indiscernibles is incorporated into set theory.

\medskip

\begin{definition}
\textnormal{(\cite[Definition 3.1]{E-ZFI}).}
The theory $\mathsf{ZFCI}_{<}$ is formulated in $\{\in,<,I\}$. In addition to
$\mathsf{ZFC}$ in this expanded language, it asserts that $<$ is a global
well-order satisfying the global-well-order axioms $\mathsf{GW}$, that $I$ is
cofinal in $\Ord$, and that $I$ is a class of order
indiscernibles for $(V,\in,<)$. Thus, for
each formula $\varphi(x_1,\ldots,x_k)$ in the language $\{\in,<\}$, and
any two increasing $k$-tuples from $I$, the truth value of $\varphi$
is the same on the two tuples.
\end{definition}

The arithmetic counterpart, $\mathsf{PAI}$, augments Peano arithmetic
with an unbounded set of order indiscernibles and induction in the expanded
language. It is conservative over $\PA$: the proof uses Ramsey's theorem
to satisfy any prescribed finite collection of indiscernibility requirements
\cite[Theorem 3.3]{E-PAI}. Via the corresponding bi-interpretability with
finite set theory, this yields the following formulation.
\begin{theorem}
For every sentence $\varphi$ of the pure language of set theory,
the following are equivalent:
\begin{enumerate}
\item $\mathsf{ZFCI}^{-\infty}_{<}\vdash\varphi$;
\item $\mathsf{ZFC}^{-\infty}\vdash\varphi$;
\item $\mathsf{GBC}^{-\infty}+$ ``$\Ord$ is weakly compact'' $\vdash\varphi$.
\end{enumerate}
\end{theorem}

In contrast, restoring Infinity changes the first-order strength of the
theory: the corresponding result leads to the Levy scheme.
\begin{theorem}
\textnormal{(\cite[Theorem 4.3]{E-ZFI}).}
For every sentence $\varphi$ of the pure language of set theory,
the following are equivalent:
\begin{enumerate}
\item $\mathsf{ZFCI}_{<}\vdash\varphi$;
\item $\mathsf{ZFC}+\Lambda\vdash\varphi$;
\item $\mathsf{GBC}+$ ``$\Ord$ is weakly compact'' $\vdash\varphi$.
\end{enumerate}
\end{theorem}

The proof that indiscernibles yield $\Lambda$ uses the
Ehrenfeucht--Mostowski stretching construction. Starting with the
Skolem hull $M_I$ of the indiscernibles, one appends a linear order $J$
without a least element, obtaining an elementary end extension
$M_{I+J}$ of $M_I$. Moreover, the argument in \cite[Lemma 4.2]{E-ZFI}
implicitly shows that, with respect to a suitable set-like global
well-ordering, $M_I$ is a strong cut of $M_{I+J}$. As in the
automorphism construction of Exhibit I, this property helps establish
that the classes coded in the larger model form a model of
$\mathsf{GBC}$ over the smaller one. The stretching construction also
yields the tree property for $\Ord$ in this class expansion, so the
Levy scheme follows by Exhibit II.

For the reverse implication, the proof in \cite{E-ZFI} uses the partition properties of $\Ord$ available in $\mathsf{GBC}+$ ``$\Ord$ is weakly compact''. By passing to a suitable elementary extension, one obtains a partition property with a nonstandard finite exponent. Together with a full satisfaction class, this yields a cofinal class of indiscernibles satisfying the required expanded-language axioms. The use of the partition properties of $\Ord$ also connects this construction with the iterated-ultrapower method underlying the automorphism results of Exhibit I.

\begin{corollary}
For every consistent completion $T$ of $\mathsf{ZFC}$, the following are
equivalent: 

\begin{enumerate}
    \item $T\vdash\Lambda$.
   \item Some model of $T$ admits
an expansion to a model of $\mathsf{ZFCI}_{<}$.
\end{enumerate}
\end{corollary}

The corollary asserts the existence of a suitable expansion of
\emph{some} model of $T$, not necessarily of every model of $T$.

\section{Exhibit IV: Quine--Jensen set theory}
\label{sec:exhibit-IV}

Our final comparison concerns a set theory which, unlike $\mathsf{ZFC}$, admits a
universal set. Quine's \emph{New Foundations} ($\mathsf{NF}$) \cite{Q} is based on extensionality
and stratified comprehension. A formula is stratified if integer types can be
assigned to its variables so that $x\in y$ requires the type of $y$ to be one
higher than that of $x$, while $x=y$ requires equal types. Thus $x=x$ is
stratified, and $\mathsf{NF}$ proves the existence of a universal set. The consistency of NF, a longstanding open problem, was recently established by Holmes and Wilshaw, with the central construction formally verified in Lean \cite{HolmesWilshaw}. Jensen's
modification $\mathsf{NFU}$ allows urelements, by restricting extensionality to sets;
its consistency was established long ago by Jensen \cite{Jen}. For a detailed account of the foundations and model theory of
$\mathsf{NFU}$, see Holmes \cite{Ho-1,Ho-2}.

The distinction between \emph{Cantorian} and \emph{strongly Cantorian} sets
is essential here. For a set $X$ write
\[
\iota``X=\{\{x\}:x\in X\}.
\]
The set $X$ is \emph{Cantorian} if there is a bijection from $X$ onto
$\iota``X$; it is \emph{strongly Cantorian} if the particular singleton
assignment $x\mapsto\{x\}$, restricted to $X$, is itself a set-function.
The distinction is meaningful because the defining formula of the singleton
assignment is not stratified. Strongly Cantorian sets
admit a controlled form of type-raising and type-lowering and thereby behave
more like sets in ordinary set theory \cite{Ho-2}. A Cantor-style diagonal argument shows that the universal set of $\mathsf{NFU}$ is not Cantorian, and a fortiori, is not strongly Cantorian.  

In the discussion below, we follow Holmes' convention and take $\mathsf{NFU}$ to include the axioms of Infinity and Choice, although some authors, following Jensen's original formulation, use $\mathsf{NFU}$ for the theory without these additional axioms.

The \emph{Axiom of Cantorian Sets}, originally proposed by C. Ward Henson,
asserts that every Cantorian set is strongly Cantorian. Following Solovay,
let
\[
\mathsf{NFUA}:=\mathsf{NFU}+
   \text{``every Cantorian set is strongly Cantorian''}.
\]
Holmes's study of stronger systems, notably $\mathsf{NFUB}$ and $\mathsf{NFUM}$, places
this axiom in a broader hierarchy of extensions of $\mathsf{NFU}$
\cite{Ho-2}; already at the level of $\mathsf{NFUA}$ the singleton
principle has considerable set-theoretic strength, as indicated in the result below.

The finite version of $\NFUA$ has a particularly close connection with
Peano arithmetic and finite set theory. In joint work with Robert Solovay
in 2002, developed through an extended email correspondence subsequently
edited and typeset by Thomas Forster \cite{ES-Corr}, we obtained the
following characterization. The main components of the proof are worked out in \cite{E-Tehran}. In part (b) of the theorem below, $CZ$ denotes
the \emph{Cantorian part of the Zermelian tower} of the ambient model (in the jargon of specialists of this subject).
\needspace{6\baselineskip}
\begin{theorem}\label{Enayat-Solovay-arith}
For a theory $T$ in the language $\{\in\}$, the following are equivalent:
\begin{enumerate}
\item $T$ is a consistent completion of $\mathsf{ZFC}^{-\infty}$.
\item There is a model $\mathcal A\models\mathsf{NFUA}^{-\infty}$ such that
$T=\operatorname{Th}(CZ^{\mathcal A})$.
\end{enumerate}
\end{theorem}

Theorem \ref{Enayat-Solovay-arith} gives far more information than an equiconsistency result: it identifies exactly
which complete set-theoretic theories arise from the Cantorian parts in
the finite setting. The corresponding characterization in the presence
of Infinity, stated below, identifies precisely the completions of $\mathsf{ZFC}+\Lambda$.

Holmes and Solovay determined the exact consistency strength of $\mathsf{NFUA}$ as follows.

\begin{theorem}
\textnormal{(Holmes and Solovay).}
$\mathsf{NFUA}$ and $\mathsf{ZFC}+\Lambda^{-}$ are equiconsistent.\footnote{The proof of the full equiconsistency remains unpublished, although it circulated among the illuminati. Holmes \cite{Ho-2} published a refinement of Solovay's argument establishing the direction from the consistency of NFUA to the consistency of $\mathsf{ZFC} + \Lambda^-$. The converse direction relies on a combinatorial theorem of Schmerl \cite{Sch-1} concerning Mahlo cardinals of finite order. Note that Theorem \ref{My_thm_NFUA} implies this equiconsistency.}
\end{theorem}

The connection with Exhibit I is supplied by automorphisms. The
Jensen--Boffa construction of models of $\mathsf{NFU}$ from automorphisms
of models of conventional set theory has a converse based on the work of
Hinnion and Solovay. For $\mathsf{NFUA}$, the correspondence can be
refined so that the Cantorian part of the associated Zermelian structure
is the fixed-point structure of a suitable automorphism. The fixed-point
characterization in Exhibit I therefore brings the Levy scheme into play.

The following theorem from \cite{E-5} characterizes the complete
first-order theories arising as theories of the Cantorian parts of models
of $\mathsf{NFUA}$. 

\begin{theorem}\label{My_thm_NFUA}
\textnormal{(\cite{E-5}).}
For a theory $T$ in the language $\{\in\}$, the following are equivalent:
\begin{enumerate}
\item $T$ is a consistent completion of $\mathsf{ZFC}+\Lambda$.
\item There is a model $\mathcal A\models \mathsf{NFUA}$ such that
$T=\operatorname{Th}(CZ^{\mathcal A})$.
\end{enumerate}
\end{theorem}

\section{A Brief History of \texorpdfstring{$\Lambda$}{Lambda}}

The two strands that meet in the Levy scheme have different origins:
the hierarchy of Mahlo cardinals and the logical study of reflection.
Mahlo's investigations, beginning in 1911 \cite{Mahlo1911},
predated the modern theory of stationary sets and the subsequent
large-cardinal hierarchy.  The use of closed unbounded and stationary
sets eventually supplied the language in which the finite iterations
of Mahloness could be formulated and studied systematically.

The reflection strand acquired a precise metamathematical form in the
work of Montague. The reflection theorem was used by Montague (together with G\"odel's second
incompleteness theorem) to show that $\mathsf{ZF}$ is not
finitely axiomatizable over Zermelo set theory; see
\cite{Montague1961,MontagueFraenkel}.

Levy then connected reflection systematically with large-cardinal
notions.  His work \cite{LevyReflection,LevyStrongInfinity} showed, in
particular, how the usual reflection scheme can be viewed as a
formulation of Infinity and Replacement over an appropriate weak
set-theoretic base, and how stronger forms of reflection lead to
inaccessible and Mahlo cardinals.  One early instance of this
connection is the equivalence, over $\mathsf{ZF}$, of the following
schemes:
\begin{enumerate}
\item For each $n<\omega$, the sentence asserting that for every ordinal
$\beta$ there is an inaccessible cardinal $\kappa>\beta$ such that
$V_\kappa\prec_{\Sigma_n}V$.
\item For each formula $\varphi(x,\vec y)$ in the language of set theory,
the sentence asserting that every club class of ordinals defined by
$\varphi(x,\vec a)$, with set parameters $\vec a$, contains an
inaccessible cardinal.
\end{enumerate}
The generalizations involving higher finite Mahlo degrees are among
the equivalent forms of $\Lambda$ discussed in Section~3.  Levy's
subsequent analysis of the $\Sigma_n$ hierarchy
\cite{LevyHierarchy} also made it possible to express each finite
level of correctness by a single first-order formula.  In the
terminology adopted here, $\Lambda$ synchronizes finite Mahloness
with finite correctness, and may suggestively be read as saying that
``$\Ord$ is $\omega$-Mahlo.''

A different impetus came from combinatorics and model theory in the
1970s.  Schmerl's analysis of partition properties of finite Mahlo
levels \cite{Sch-1} was subsequently developed further, notably by
Hajnal, Kanamori, and Shelah \cite{HKS}.  In parallel, Schmerl and
Shelah \cite{SS} investigated the Hanf-number phenomenon for
power-like models.  Their results showed how finite levels of
Mahloness govern the transfer from the existence of sufficiently
large power-like models to the existence of models at other
uncountable cardinalities.

The investigation of elementary end extensions supplied another
closely related route.  Kaufmann \cite{Kau} obtained extension and
obstruction results that sharpened the contrast with the
MacDowell--Specker theorem for arithmetic.  Kaufmann's paper was a source of inspiration for me, especially since it was there that I first encountered the Levy scheme, in a proof outline of Theorem~\ref{Lambda-eee}, under the additional assumption that Mahlo cardinals of all finite orders exist in the sense of the metatheory. This additional assumption was eventually removed in \cite{E-5}. My subsequent work
\cite{E-1,E-2,E-3,E-4,E-5,{E-ZFI}} brought the Levy scheme into increasingly
precise characterizations involving conservative extensions,
power-like models, and automorphisms.

Another route from the model theory of
arithmetic into set theory was developed by McAloon and Ressayre
\cite{McAloonRessayre}, who adapted the Kirby--Paris method of indicators
to models of set theory.  Their work connected combinatorial partition
principles and finite levels of the Mahlo hierarchy with special
initial segments and extensions of models of set theory, thereby
transporting to set theory methods that had proved effective in the
study of cuts and independence phenomena in models of arithmetic.
This line is especially relevant here because it also leads to
conservative elementary rank extensions under class-theoretic
partition hypotheses.  More recently, De Bondt \cite{DeBondt} gave a detailed
reconstruction of one such argument from a proof sketch of McAloon
and Ressayre: for a countable model of $\mathsf{GBC}$ satisfying the
partition property $\Ord\longrightarrow(\Ord)^3_2$, the associated
structure expanded by predicates for its classes has a conservative
elementary rank extension.
 
As pointed out in the introduction, the finite Mahlo hierarchy has an independent source of
motivation in Friedman's investigations of concrete independence
phenomena. However, to my knowledge there is no direct link between the machinery developed by Friedman and those surveyed in this paper.


\begin{thebibliography}{99}

\bibitem{Ackermann}W. Ackermann, \textit{Die Widerspruchsfreiheit der allgemeinen Mengenlehre}, \textbf{Mathematische Annalen} 114 (1937), 305--315. DOI: 10.1007/BF01594179.

\bibitem{Bo}W. Boos, \textit{Boolean extensions which efface the Mahlo property}, \textbf{Journal of Symbolic Logic} 39 (1974), no. 2, 254--268. DOI: 10.2307/2272638.

\bibitem{JensenSolovay}R. B. Jensen and R. M. Solovay, \textit{Some applications of almost disjoint sets}, in \textit{Mathematical Logic and Foundations of Set Theory} (Proc. Internat. Colloq., Jerusalem, 1968), North-Holland, Amsterdam, 1970, pp. 84--104.

\bibitem{DeBondt}B. De Bondt, \textit{Sets, Trees and Incompleteness}, Master's thesis, Faculty of Sciences, Ghent University, 2018.

\bibitem{E-1}A. Enayat, \textit{On certain elementary extensions of models of set theory}, \textbf{Transactions of the American Mathematical Society} 283 (1984), 705--715.

\bibitem{E-2}A. Enayat, \textit{Conservative elementary extensions of models of set theory and generalizations}, \textbf{Journal of Symbolic Logic} 51 (1986), 1005--1021.

\bibitem{E-3}A. Enayat, \textit{Analogues of the MacDowell--Specker theorem in set theory}, in \textit{Models, Algebras and Proofs}, edited by X. Caicedo and C. H. Montenegro, Marcel Dekker, 1998, pp. 25--50.

\bibitem{E-4}A. Enayat, \textit{Power-like models of set theory}, \textbf{Journal of Symbolic Logic} 66 (2001), no. 4, 1766--1782. DOI: 10.2307/2694973.

\bibitem{E-5}A. Enayat, \textit{Automorphisms, Mahlo cardinals, and NFU}, in \textit{Nonstandard Models of Arithmetic and Set Theory}, edited by A. Enayat and R. Kossak, Contemporary Mathematics 361, American Mathematical Society, 2004, pp. 37--59.

\bibitem{E-Tehran}
A. Enayat, \textit{From bounded arithmetic to second order
arithmetic via automorphisms},
in \textit{Logic in Tehran},
edited by A. Enayat, I. Kalantari, and M. Moniri,
Lecture Notes in Logic 26,
Association for Symbolic Logic and A K Peters,
2006, pp. 87--113.
DOI: 10.1017/9781316755747.008.

\bibitem{E-2009}A. Enayat, \textit{An inevitable extension of ZFC}, lecture presented at Kunen Fest, University of Wisconsin--Madison, April 3, 2009. DOI: 10.13140/RG.2.1.4897.0324.

\bibitem{E-Visser}A. Enayat, \textit{Variations on a Visserian theme}, in \textit{Liber Amicorum Alberti: A Tribute to Albert Visser}, edited by J. van Eijck, R. Iemhoff, and J. J. Joosten, College Publications, London, 2016, pp. 99--110.

\bibitem{E-ZFI}A. Enayat, \textit{Set theory with a proper class of indiscernibles}, \textbf{Fundamenta Mathematicae} 259 (2022), 33--76. DOI: 10.4064/fm999-2-2022.

\bibitem{E-PAI}A. Enayat, \textit{Indiscernibles and satisfaction classes in arithmetic}, \textbf{Archive for Mathematical Logic} 63 (2024), 655--677. DOI: 10.1007/s00153-024-00915-8.

\bibitem{E-Visser-Corr}A. Enayat, \textit{Corrigendum \& Addendum to: Variations on a Visserian theme}, arXiv:2602.15081 (2026).

\bibitem{E-Deadends}A. Enayat, \textit{Models of set theory: Extensions and dead-ends}, \textbf{Journal of Symbolic Logic}, First View (published online 24 March 2026), 1--38. DOI: 10.1017/jsl.2026.10199.

\bibitem{ES-Corr}
A. Enayat and R. M. Solovay,
\textit{Solovay--Enayat correspondence from 2002},
unpublished correspondence, edited and typeset by T. E. Forster,
January 8, 2022, 64 pp.

\bibitem{EH-Ord}A. Enayat and J. D. Hamkins, \textit{ZFC proves that the class of ordinals is not weakly compact for definable classes}, \textbf{Journal of Symbolic Logic} 83 (2018), no. 1, 146--164. DOI: 10.1017/jsl.2017.75.

\bibitem{ESV}A. Enayat, J. H. Schmerl, and A. Visser, \textit{$\omega$-models of finite set theory}, in \textit{Set Theory, Arithmetic, and Foundations of Mathematics: Theorems, Philosophies}, edited by J. Kennedy and R. Kossak, Lecture Notes in Logic 36, Cambridge University Press, 2011, pp. 43--65.

\bibitem{Friedman1981}H. M. Friedman, \textit{On the necessary use of abstract set theory}, \textbf{Advances in Mathematics} 41 (1981), no. 3, 209--280. DOI: 10.1016/0001-8708(81)90021-9.

\bibitem{Friedman2004}H. M. Friedman, \textit{Working with nonstandard models}, in \textit{Nonstandard Models of Arithmetic and Set Theory}, A. Enayat and R. Kossak (eds.), Contemporary Mathematics 361, American Mathematical Society, Providence, RI, 2004, pp. 71--86. DOI: 10.1090/conm/361/06589. MR 2089277.

\bibitem{G}H. Gaifman, \textit{Models and types of arithmetic}, \textbf{Annals of Mathematical Logic} 9 (1976), 223--306.

\bibitem{GodelContinuum}K. G\"odel, \textit{What is Cantor's continuum problem?}, \textbf{American Mathematical Monthly} 54 (1947), 515--525; revised and expanded version in \textit{Philosophy of Mathematics: Selected Readings}, edited by P. Benacerraf and H. Putnam, Prentice--Hall, 1964.

\bibitem{HKS}A. Hajnal, A. Kanamori, and S. Shelah, \textit{Regressive partition relations for infinite cardinals}, \textbf{Transactions of the American Mathematical Society} 299 (1987), no. 1, 145--154. DOI: 10.1090/S0002-9947-1987-0869404-5.

\bibitem{Ho-1}R. Holmes, \textit{Elementary Set Theory with a Universal Set}, Cahiers du Centre de Logic, vol. 10, Academia, Louvain-la-Neuve, Belgium, 1998.

\bibitem{Ho-2}M. Randall Holmes, \textit{Strong axioms of infinity in NFU}, \textbf{Journal of Symbolic Logic} 66 (2001), no. 1, 87--116. DOI: 10.2307/2694912. See also the erratum, \textbf{Journal of Symbolic Logic} 66 (2001), no. 4, p. 1974.

\bibitem{HolmesWilshaw}
M. R. Holmes and S. Wilshaw,
\textit{NF is Consistent},
arXiv:1503.01406 [math.LO], revised 2025.

\bibitem{Jec}T. Jech, \textit{Set Theory}, Academic Press, New York, 1978.

\bibitem{Jen}R. B. Jensen, \textit{On the consistency of a slight modification of Quine's New Foundations}, \textbf{Synthese} 19 (1969), 250--263.

\bibitem{Kan}A. Kanamori, \textit{The Higher Infinite}, Springer-Verlag, Berlin, 1994.

\bibitem{Kau}M. Kaufmann, \textit{Blunt and topless end extensions of models of set theory}, \textbf{Journal of Symbolic Logic} 48 (1983), 1053--1073.

\bibitem{KW}R. Kaye and T. L. Wong, \textit{On interpretations of arithmetic and set theory}, \textbf{Notre Dame Journal of Formal Logic} 48 (2007), 497--510.

\bibitem{Keisler-Morley}H. J. Keisler and M. Morley, \textit{Elementary extensions of models of set theory}, \textbf{Israel Journal of Mathematics} 5 (1968), 49--65.

\bibitem{Keisler-Silver}H. J. Keisler and J. Silver, \textit{End extensions of models of set theory}, in \textit{Axiomatic Set Theory (Part I)}, Proceedings of Symposia in Pure Mathematics 13, American Mathematical Society, 1970, pp. 177--187.

\bibitem{LevyReflection}A. Levy, \textit{Principles of reflection in axiomatic set theory}, \textbf{Fundamenta Mathematicae} 49 (1960), 1--10.

\bibitem{LevyStrongInfinity}A. Levy, \textit{Axiom schemata of strong infinity in axiomatic set theory}, \textbf{Pacific Journal of Mathematics} 10 (1960), 223--238. DOI: 10.2140/pjm.1960.10.223.
\bibitem{LevyHierarchy}A. Levy, \textit{A hierarchy of formulas in set theory}, Memoirs of the American Mathematical Society 57, 1965.




\bibitem{LevySolovay}A. Levy and R. M. Solovay, \textit{Measurable cardinals and the continuum hypothesis}, \textbf{Israel Journal of Mathematics} 5 (1967), no. 4, 234--248. DOI: 10.1007/BF02771612.

\bibitem{MacDowell-Specker}R. MacDowell and E. Specker, \textit{Modelle der Arithmetik}, in \textit{Infinitistic Methods}, Pergamon Press, 1961, pp. 257--263.

\bibitem{Mahlo1911}P. Mahlo, \textit{\"Uber lineare transfinite Mengen}, \textbf{Berichte \"uber die Verhandlungen der K\"oniglich S\"achsischen Gesellschaft der Wissenschaften zu Leipzig, Math.-Phys. Klasse} 63 (1911), 187--225.

\bibitem{McAloonRessayre}K. McAloon and J.-P. Ressayre, \textit{Les m\'ethodes de Kirby--Paris et la th\'eorie des ensembles}, in C. Berline, K. McAloon, and J.-P. Ressayre (eds.), \textit{Model Theory and Arithmetic}, Lecture Notes in Mathematics 890, Springer-Verlag, Berlin, 1981, pp. 154--184.

\bibitem{Montague1961}R. Montague, \textit{Semantical closure and non-finite axiomatizability I}, in \textit{Infinitistic Methods}, Pergamon Press, 1961, pp. 45--69.

\bibitem{MontagueFraenkel}R. Montague, \textit{Fraenkel's addition to the axioms of Zermelo}, in \textit{Essays on the Foundations of Mathematics}, Magnes Press, Jerusalem, 1961, pp. 91--114.

\bibitem{Mycielski}J. Mycielski, \textit{The definition of arithmetic operations in the Ackermann model}, \textbf{Algebra i Logika Seminar} 3 (1964), 64--65. [In Russian.]

\bibitem{Phillips}R. G. Phillips, \textit{Omitting types in arithmetic and conservative extensions}, in \textit{Victoria Symposium on Nonstandard Analysis}, Lecture Notes in Mathematics 369, Springer, 1974, pp. 195--202.

\bibitem{Q}W. V. O. Quine, \textit{New foundations for mathematical logic}, \textbf{American Mathematical Monthly} 44 (1937), 70--80.

\bibitem{Sch-1}J. H. Schmerl, \textit{A partition property characterizing cardinals hyperinaccessible of finite type}, \textbf{Transactions of the American Mathematical Society} 188 (1974), 281--291.

\bibitem{SS}J. Schmerl and S. Shelah, \textit{On power-like models for hyperinaccessible cardinals}, \textbf{Journal of Symbolic Logic} 37 (1972), 531--537.

\end{thebibliography}
\end{document}